# OBSTRUCTION THEORY IN MODEL CATEGORIES

J. DANIEL CHRISTENSEN, WILLIAM G. DWYER, AND DANIEL C. ISAKSEN

ABSTRACT. Many examples of obstruction theory can be formulated as the study of when a lift exists in a commutative square. Typically, one of the maps is a cofibration of some sort and the opposite map is a fibration, and there is a functorial obstruction class that determines whether a lift exists. Working in an arbitrary pointed proper model category, we classify the cofibrations that have such an obstruction theory with respect to all fibrations. Up to weak equivalence, retract, and cobase change, they are the cofibrations with weakly contractible target. Equivalently, they are the retracts of principal cofibrations. Without properness, the same classification holds for cofibrations with cofibrant source. Our results dualize to give a classification of fibrations that have an obstruction theory.

## 1. INTRODUCTION

The following extension-lifting problem is ubiquitous in modern homotopy theory. Consider a commutative square

$$\begin{array}{ccc} A & \longrightarrow & X \\ i \downarrow & & \downarrow p \\ B & \longrightarrow & Y \end{array}$$

in which $i$ is a cofibration (or less technically, some kind of monomorphism) while $p$ is a fibration (or some kind of epimorphism). When does a map $B \to X$ exist making both triangles commute?

Classical obstruction theory [14] [15] gives a detection principle for existence (and uniqueness) of lifts in the category of spaces in terms of homotopy theory.

In this paper, we show how some aspects of classical obstruction theory are entirely abstract, working the same for any model category. In addition to providing tools for lifting in non-topological contexts, this approach enlightens classical obstruction theory by showing that much of it does not depend on specific properties of topological spaces.

We start with a pointed model category $\mathcal{C}$. Examples include pointed topological spaces, pointed simplicial sets, spectra, and chain complexes of modules over a ring. Fix a cofibration $i$. We say that $i$ has an obstruction theory if there exists some object $W$ in $\mathcal{C}$ such that for every commutative square as above, there is a well-defined weak homotopy class from $W$ to the fibre of $p$ (called the obstruction), such

Date: November 4, 2018.
1991 *Mathematics Subject Classification.* 55S35, 55U35, 18G55 (primary); 18G30, 55P42 (secondary).
The first author was supported by an NSERC grant.
The third author was supported by an NSF Postdoctoral Research Fellowship.

that a lift exists in the square if and only if the map from $W$ to the fibre is weakly null-homotopic. We also require that the obstruction class be functorial in the fibration $p$ (see Definition 4.1).

Starting from this very general framework, we prove a theorem precisely classifying the cofibrations that have obstruction theories. Assuming that $\mathcal{C}$ is both left and right proper, the class of cofibrations having obstruction theories is the smallest class of cofibrations containing all cofibrations with weakly contractible target and closed under retract, weak equivalence, and cobase change. In other words, they are the retracts of principal cofibrations. Although $W$ is not uniquely determined by the definitions, when $i$ is a principal cofibration, $W$ can be taken to be a desuspension of the cofibre (see Remark 7.3). Moreover, the category of obstruction theories is contractible (see Subsection 4.2).

Our theorem explains why some kind of hypotheses are needed in classical obstruction theory. It is not possible to give an obstruction theory for topological spaces without some kind of restriction to principal cofibrations or fibrations.

We obtain slightly weaker results when $\mathcal{C}$ is not proper. Without left properness, we must assume that the source of the cofibration is cofibrant. Without right properness, we must assume that the targets of the fibrations are fibrant. Here the notion of a fibrant obstruction theory becomes useful. The essential reason that fibrancy and cofibrancy assumptions take the place of properness is that base changes along fibrations preserve weak equivalences between fibrant objects and cobase changes along cofibrations preserve weak equivalences between cofibrant objects. The point is that some result about base changes and cobase changes of weak equivalences is necessary; this result either comes by assumption from properness or from having enough cofibrancy and fibrancy.

Obstruction theories in the sense of this paper have been used for simplicial sets with its usual model structure in [6] and [10]. The obstruction for lifting a standard generating cofibration with respect to a fibration is an element of a homotopy group of a fibre. See Section 8.1 for more details. This paper grew out of understanding precisely how this special case works.

In a stable model category (see Section 8.3) we show that every cofibration has an obstruction theory. In the special case of the model category of unbounded chain complexes of objects in an abelian category, the results of the present paper provide a conceptual explanation for results such as [3, Lemma 2.3]. For spectra, obstruction theory for the standard generating cofibrations is used in [12]. See Section 8.3 for more details. In addition, in the stable model category of pro-spectra, the results of this paper produce new tools for handling lifting problems, tools which were used in early versions of [4].

We work exclusively on the question of when cofibrations have obstruction theories, but everything dualizes to the study of fibrations that have co-obstruction theories. A fibration $p$ has a co-obstruction theory if there exists some object $W$ in $\mathcal{C}$ such that for every commutative square as above, there is a well-defined weak homotopy class from the cofibre of $i$ to $W$ (called the obstruction), such that a lift exists in the square if and only if the map from the cofibre to $W$ is weakly null-homotopic. We also require that the obstruction class be functorial in the cofibration $i$. Co-obstruction theory is a key technique of [11]. Also, this notion of co-obstruction theory is precisely what appears in [14, § 8.2] for the special case of a principal fibration of topological spaces whose fibre is an Eilenberg-Mac Lane space $K(\pi, n-1)$. In this case, $W$ is the delooping $K(\pi, n)$ of the fibre, and the obstruction lies in $[B/A, K(\pi, n)] = H^n(B, A; \pi)$, where the cofibration is $A \to B$. We do not discuss the iterative procedure that occurs when lifting against a tower of fibrations.

We are curious how many other examples of obstruction theory can be viewed as special cases of our theory, but have not yet investigated this in detail.

A summary of the contents of the paper follows. We begin in Section 2 by introducing the category of arrows $\text{Ar}\mathcal{C}$ of a pointed model category $\mathcal{C}$ and describing two useful model structures on $\text{Ar}\mathcal{C}$. Then in Section 3 we establish some technical lifting results. In Section 4 we give the definitions of obstruction theory, fibrant obstruction theory, and cofibrant fibrant obstruction theory and describe their elementary properties. Subsections 4.1 and 4.2 discuss a rigidification of the notion of obstruction theory and the uniqueness of obstruction theories for a given cofibration.

In Sections 5 and 6 we exhibit several ways of producing cofibrations that have obstruction theories. First, cofibrations with weakly contractible target have obstruction theories. Second, retracts and cobase changes preserve cofibrations that have obstruction theories. Third, and most difficult, weak equivalences preserve cofibrations that have obstruction theories. This allows us to prove the main classification theorem in Section 7.

Section 8 contains some applications. We explain how to apply our notions of obstruction theory to the unpointed category of simplicial sets or of topological spaces. We also show that every cofibration in a stable model category has an obstruction theory.

We assume that the reader is familiar with model categories. The original reference is [13], but see [7], [8], or [5] for more modern treatments. We follow the notation and conventions of [7] as closely as possible.

## 2. The Category of Morphisms

Throughout the entire paper, $\mathcal{C}$ denotes a pointed model category with functorial factorizations. Let $\text{Ar}\mathcal{C}$ be the category of morphisms in $\mathcal{C}$. It is a diagram category, where the index category has two objects and a unique map between them. Thus objects of $\text{Ar}\mathcal{C}$ are morphisms in $\mathcal{C}$, and morphisms in $\text{Ar}\mathcal{C}$ are commuting squares in $\mathcal{C}$. When the meaning is clear, we write $X$ for the identity object $X \to X$ of $\text{Ar}\mathcal{C}$.

The category $\text{Ar}\mathcal{C}$ supports two model structures: the injective structure and the projective structure [2, p. 314]. In both structures, the weak equivalences are levelwise weak equivalences. Therefore, the associated homotopy categories are identical.

**Definition 2.1.** A **weak equivalence** (or more precisely a **levelwise weak equivalence**) (*resp.*, **injective cofibration**, **projective fibration**) from a map $f: X \to X'$ to another map $g: Y \to Y'$ in the category $\text{Ar}\mathcal{C}$ is a commutative square

$$\begin{array}{ccc} X & \longrightarrow & Y \\ f \downarrow & & \downarrow g \\ X' & \longrightarrow & Y' \end{array}$$

in which the horizontal maps are weak equivalences (*resp.*, cofibrations, fibrations). The map $f \to g$ is a **projective cofibration** if both $X \to Y$ and $X' \amalg_X Y \to Y'$ are cofibrations. The map $f \to g$ is an **injective fibration** if both $X' \to Y'$ and $X \to X' \times_{Y'} Y$ are fibrations.

Both model structures are examples of Reedy model structures [7, Ch. 17]. One way of obtaining the injective structure is by considering the opposite of the projective model structure on $\mathrm{Ar}(\mathcal{C}^{\mathrm{op}})$.

Note that all projective cofibrations are injective cofibrations. Similarly, all injective fibrations are projective fibrations. A cofibration in $\mathcal{C}$ is an injective cofibrant object if and only if it is a projective cofibrant object if and only if it has cofibrant source. Similarly, a fibration in $\mathcal{C}$ is projective fibrant if and only if it is injective fibrant if and only if it has fibrant target.

We restate here the following useful lemma about base changes and cobase changes of weak equivalences in model categories.

**Lemma 2.2.** *Base changes along fibrations preserve weak equivalences between fibrant objects. Dually, cobase changes along cofibrations preserve weak equivalences between cofibrant objects.*

*Proof.* See [7, Prop. 11.1.2]. □

The following lemma is used many times throughout the paper. It allows us to unify arguments that either assume fibrancy or right properness.

**Lemma 2.3.** *Let $p \to p'$ be a weak equivalence between fibrations $p : X \to Y$ and $p' : X' \to Y'$ with fibres $F$ and $F'$ respectively. If $Y$ and $Y'$ are fibrant, or $\mathcal{C}$ is right proper, then the induced map $F \to F'$ is a weak equivalence.*

*Proof.* Let $P$ be the fibre product $X' \times_{Y'} Y$. Consider the diagram

$$\begin{array}{ccccc} X & \longrightarrow & P & \longrightarrow & X' \\ p \downarrow & & \downarrow & & \downarrow p' \\ Y & \longrightarrow & Y & \longrightarrow & Y'. \end{array}$$

The map $P \to X'$ is a weak equivalence: in the right proper case, this is by definition, and in the case where $Y$ and $Y'$ are fibrant, this follows from Lemma 2.2. Therefore the map $X \to P$ is also a weak equivalence by the two-out-of-three axiom. The fibre of $P \to Y$ is isomorphic to $F'$. Hence it suffices to assume that $Y$ equals $Y'$.

Consider the model category $\mathcal{C} \downarrow Y$ of objects over $Y$. Let $R : \mathcal{C} \downarrow Y \to \mathcal{C}$ be the functor taking a map $A \to Y$ to $* \times_Y A$. It is straightforward to check that $R$ is right adjoint to the functor $L : \mathcal{C} \to \mathcal{C} \downarrow Y$ taking an object $A$ to the trivial map $A \to Y$.

The functor $L$ preserves cofibrations and acyclic cofibrations, so $L$ and $R$ form a Quillen pair. This means that $R$ preserves weak equivalences between fibrant objects. The fibrant objects of $\mathcal{C} \downarrow Y$ are the fibrations with target $Y$, so $p$ and $p'$ are both fibrant in this category. Therefore, $Rp \to Rp'$ is a weak equivalence.

The right proper case also follows from [7, Prop. 11.2.9]. □

Frequently in abstract homotopy theory, one has a square

$$\begin{array}{ccc} A & \longrightarrow & X \\ i \downarrow & & \downarrow p \\ B & \longrightarrow & Y \end{array}$$

commuting in $\mathcal{C}$ and wants to know whether a map $B \to X$ exists making both resulting triangles commute. We rewrite this problem in terms of the category $\mathrm{Ar}\mathcal{C}$.

**Lemma 2.4.** *Suppose given a square as above in the category* $\mathcal{C}$. *A lift exists for this square if and only if a lift exists in the square*

$$\begin{array}{ccc} * & \longrightarrow & X \\ \downarrow & & \downarrow \\ i & \longrightarrow & p \end{array}$$

*in* $\operatorname{Ar}\mathcal{C}$ *if and only if a lift exists in the square*

$$\begin{array}{ccc} i & \longrightarrow & p \\ \downarrow & & \downarrow \\ B & \longrightarrow & * . \end{array}$$

*Proof.* The proof is a straightforward diagram chase. □

We shall frequently switch between these three equivalent forms of the lifting problem.

## 3. Lifting Results

We now study how lifts for a given cofibration carry over to another weakly equivalent cofibration.

**Proposition 3.1.** *Let* $i : A \to B$ *and* $i' : A' \to B'$ *be cofibrations, and let* $i \to i'$ *be a weak equivalence in* $\operatorname{Ar}\mathcal{C}$. *Also assume that* $A$ *and* $A'$ *are cofibrant or that* $\mathcal{C}$ *is left proper. Let* $p : X \to Y$ *be any fibration. Then a lift exists in the square*

(3.1)
$$\begin{array}{ccc} A' & \longrightarrow & X \\ i' \downarrow & & \downarrow p \\ B' & \longrightarrow & Y \end{array}$$

*if and only if a lift exists in the square*

(3.2)
$$\begin{array}{ccccc} A & \xrightarrow{\sim} & A' & \longrightarrow & X \\ i \downarrow & & & & \downarrow p \\ B & \xrightarrow{\sim} & B' & \longrightarrow & Y. \end{array}$$

*Proof.* The proof is very similar to the proof of [7, Prop. 11.1.16]. If a lift $B' \to X$ exists in the square 3.1, then the composition $B \to B' \to X$ is a lift for the square 3.2. The other implication is harder.

Suppose that a lift $B \to X$ exists for the square 3.2. Let $P$ be the pushout $A \amalg_B A'$. Because of left properness or because of Lemma 2.2, the map $B \to P$ is a weak equivalence, so $P \to B'$ is a weak equivalence by the two-out-of-three axiom. Since $j : A' \to P$ is a cobase change of $i$, there is a lift $P \to X$ in the square

$$\begin{array}{ccc} A' & \longrightarrow & X \\ j \downarrow & & \downarrow p \\ P & \longrightarrow & Y. \end{array}$$

Now consider the model category $A' \downarrow \mathcal{C} \downarrow Y$ of objects under $A'$ and over $Y$. Then $P$ and $B'$ (equipped with the obvious structure maps from $A'$ and to $Y$) are cofibrant objects of this category, while $X$ is a fibrant object. Moreover, there are morphisms $P \to B'$ and $P \to X$ in $A' \downarrow \mathcal{C} \downarrow Y$. The first map is a weak equivalence, so there is also a map $B' \to X$ in $A' \downarrow \mathcal{C} \downarrow Y$ by [7, Cor. 7.6.5]. □

The point of this result is that finding lifts for $i'$ reduces to finding lifts for $i$. Given any lifting problem for $i'$ (in other words, a square 3.1), we can consider the square 3.2 instead. However, finding lifts for $i$ does not reduce to finding lifts for $i'$: in general, not every square

$$\begin{array}{ccc} A & \longrightarrow & X \\ i \downarrow & & \downarrow p \\ B & \longrightarrow & Y \end{array}$$

can be rewritten in the form 3.2, even if $A$ is cofibrant.

For future reference, we record the dual result here.

**Proposition 3.2.** *Let $p : X \to Y$ and $p' : X' \to Y'$ be fibrations, and let $p \to p'$ be a weak equivalence in $\mathrm{Ar}\mathcal{C}$. Also assume that $Y$ and $Y'$ are fibrant or that $\mathcal{C}$ is right proper. Let $i : A \to B$ be any cofibration. Then a lift exists in the square*

(3.3)
$$\begin{array}{ccc} A & \longrightarrow & X \\ i \downarrow & & \downarrow p \\ B & \longrightarrow & Y \end{array}$$

*if and only if a lift exists in the square*

(3.4)
$$\begin{array}{ccccc} A & \longrightarrow & X & \xrightarrow{\sim} & X' \\ i \downarrow & & & & \downarrow p' \\ B & \longrightarrow & Y & \xrightarrow{\sim} & Y'. \end{array}$$

*Proof.* The proof is dual to the proof of Proposition 3.1. □

## 4. Obstruction Theory

**Definition 4.1.** A cofibration $i$ **has an obstruction theory** if there exists an object $W$ such that for every fibration $p$ with fibre $F$ and every map $i \to p$ in $\mathrm{Ar}\mathcal{C}$ there exists a well-defined **obstruction** $\alpha$, which is a weak homotopy class from $W$ to $F$ (*i.e.*, an element of $[W, F]$), with the following two properties. First, a lift exists in the square

$$\begin{array}{ccc} A & \longrightarrow & X \\ i \downarrow & & \downarrow p \\ B & \longrightarrow & Y \end{array}$$

if and only if $\alpha$ is the trivial homotopy class. Second, $\alpha$ is functorial in the following sense. Given two fibrations $p$ and $p'$ with fibres $F$ and $F'$ respectively and a map $p \to p'$, the obstruction $\alpha'$ of the composition $i \to p \to p'$ is the composition of the obstruction $\alpha$ of the map $i \to p$ with the map $F \to F'$.

*Remark* 4.2. The functoriality of obstructions can be reexpressed in terms of the Grothendieck construction of the functor

$$\mathcal{C} \to \mathrm{Sets} : X \mapsto [W, X].$$

**Definition 4.3.** A cofibration **has a fibrant obstruction theory** if the conditions of Definition 4.1 (including functoriality) are satisfied for all maps $i \to p$ for which $p$ is a fibration with fibrant target. A cofibration **has a cofibrant fibrant obstruction theory** if the conditions of Definition 4.1 (including functoriality) are satisfied for all injective cofibrations $i \to p$ for which $p$ is a fibration with fibrant target.

As in Remark 4.2, the functoriality of fibrant obstruction theories and cofibrant fibrant obstruction theories can be expressed in terms of Grothendieck constructions.

These three notions of obstruction theories have some obvious relationships. If a cofibration has an obstruction theory, then it necessarily has a fibrant obstruction theory. Also, if it has a fibrant obstruction theory, then it has a cofibrant fibrant obstruction theory. Next we state some less obvious connections.

**Proposition 4.4.** *If $\mathcal{C}$ is right proper, then a cofibration has a fibrant obstruction theory if and only if it has an obstruction theory.*

*Proof.* The axioms for an obstruction theory are stronger than the axioms for a fibrant obstruction theory, so one implication follows from the definitions.

Suppose that a cofibration $i$ has a fibrant obstruction theory. Now consider a map $i \to p$ in which $p$ is a fibration between not necessarily fibrant objects. Let $\hat{p}$ be an injective fibrant replacement for $p$, so $\hat{p}$ is a fibration between fibrant objects and there is a weak equivalence $p \to \hat{p}$ in Ar$\mathcal{C}$. Let $F$ and $\hat{F}$ be the fibres of $p$ and $\hat{p}$ respectively, and let $\hat{\alpha}$ be the obstruction for the composition $i \to p \to \hat{p}$. By Lemma 2.3, the map $F \to \hat{F}$ is a weak equivalence. Define the obstruction $\alpha$ for $i \to p$ to be composition of $\hat{\alpha}$ with the weak homotopy inverse of the map $F \to \hat{F}$. This definition of $\alpha$ is functorial because injective fibrant replacements are functorial.

By Proposition 3.2, a lift exists in the square $i \to p$ if and only if a lift exists in the square $i \to p \to \hat{p}$; here we use that $\mathcal{C}$ is right proper. A lift exists in the square $i \to p \to \hat{p}$ if and only if $\hat{\alpha}$ is trivial. Since $F \to \hat{F}$ is a weak equivalence, $\hat{\alpha}$ is trivial if and only if $\alpha$ is trivial. □

The next proposition tells us that fibrant obstruction theories and cofibrant fibrant obstruction theories are actually equivalent.

**Proposition 4.5.** *A cofibration has a fibrant obstruction theory if and only if it has a cofibrant fibrant obstruction theory.*

*Proof.* The axioms for a fibrant obstruction theory are stronger than the axioms for a cofibrant fibrant obstruction theory, so one implication follows from the definitions.

Suppose that a cofibration $i$ has a cofibrant fibrant obstruction theory. For every map $i \to p$ for which $p$ is a fibration with fibrant target, take a functorial factorization

$$i \to p' \to p$$

where the first map is an injective cofibration and the second map is an acyclic injective fibration. Then $p'$ is a fibration with fibrant target, so there exists an obstruction $\alpha'$ for $i \to p'$.

Let $F'$ and $F$ be the fibres of $p'$ and $p$ respectively. By Lemma 2.3, the map $F' \to F$ is a weak equivalence; here we use that $p'$ and $p$ have fibrant targets. Define the obstruction $\alpha$ for $i \to p$ to be the composition of $\alpha'$ with the map $F' \to F$. This definition is functorial since factorizations in the injective model structure on Ar$\mathcal{C}$ are functorial.

By Proposition 3.2, a lift exists for $i \to p$ if and only if a lift exists for $i \to p'$. By assumption, a lift exists for $i \to p'$ if and only if $\alpha'$ is trivial. Since $F' \to F$ is a weak equivalence, $\alpha'$ is trivial if and only if $\alpha$ is trivial. □

The following proposition establishes a certain homotopy invariance property of obstructions. The result is not needed later, but we include it for completeness.

**Proposition 4.6.** *Let $i$ be a cofibration that has a fibrant obstruction theory, and let $p$ be a fibration with fibrant target. Suppose given two maps $f$ and $g$ (i.e., commuting squares) from $i$ to $p$ that are right homotopic in the injective or projective model structure. Then the obstructions for $f$ and for $g$ are equal.*

*Proof.* Let $p^I$ be a good path object for $p$ in the injective model structure. This is also a good path object for $p$ in the projective model structure. So in either case there exists a good right homotopy $i \to p^I$ between $f$ and $g$ [5, § 4.12]. This means that there is a commutative diagram

$$\begin{array}{c} \phantom{x} \\ i \xrightarrow{\phantom{xx}} p^I \end{array}$$

with $f$ going to $p$ above and $g$ going to $p$ below.

Note that $p^I$ is a fibration with fibrant target because $p^I \to p \times p$ is an injective fibration. There is an obstruction $\alpha$ for the map $i \to p^I$, which is a homotopy class into the fibre $F^I$ of $p^I$.

Let $F$ be the fibre of $p$. By Lemma 2.3, the map $F \to F^I$ is a weak equivalence since $p \to p^I$ is a weak equivalence. Therefore, both maps $F^I \to F$ are equal in the homotopy category since they have the same right inverse.

Since the fibrant obstruction theory for $i$ is functorial, the obstructions for $f$ and $g$ are the compositions of $\alpha$ with the two maps $F^I \to F$. Therefore, the obstructions for $f$ and $g$ are equal. □

4.1. **Rigid Obstruction Theories.** For any map $f : X \to Y$, we can choose a functorial factorization of $f$ into an acyclic cofibration $X \to X'$ followed by a fibration $f' : X' \to Y$. Write hofib for the functor $\text{Ar}\mathcal{C} \to \mathcal{C}$ sending $f$ to the fibre of $f'$. When $f$ is a fibration, there is a natural weak equivalence $\text{fib}(f) \to \text{hofib}(f)$, where $\text{fib}(f)$ denotes the fibre of $f$ (see the proof of Lemma 2.3).

Now consider a cofibration $i$ that has an obstruction theory. Then we have an obstruction $\alpha$ in $[W, \text{hofib}(i)]$ for the square $i \to i'$. Any other square $i \to p$, where $p$ is a fibration, factors through the square $i \to i'$. Therefore, the class $\alpha$ determines by functoriality the entire obstruction theory for $i$.

This implies that any obstruction theory can be rigidified in the following way. Assume without loss of generality that $W$ is cofibrant, and fix a representative $a : W \to \text{hofib}(i)$ for the class $\alpha$. To each fibration $p$ and square $i \to p$, assign the map $W \to \text{hofib}(p)$ which is the composite of $a$ with the map $\text{hofib}(i) \to \text{hofib}(p)$. This is functorial at the model category level, and via the natural weak equivalences $\text{fib}(p) \to \text{hofib}(p)$ specializes to the obstruction theory we started with.

These ideas lead to the following definition. A **rigid obstruction theory** for a cofibration $i$ is a cofibrant object $W$ and a map $a : W \to \text{hofib}(i)$ such that for each square $i \to p$ with $p$ a fibration, the square has a lift if and only if the composite $W \to \text{hofib}(i) \to \text{hofib}(p)$ is null-homotopic. A map of strict obstruction theories for $i$ from $a : W \to \text{hofib}(i)$ to $a' : W' \to \text{hofib}(i)$ is simply a map $W \to W'$ over $\text{hofib}(i)$.

4.2. **Uniqueness.** Our definition of obstruction theory does not uniquely determine $W$ nor the assignment of obstruction classes to each square.

Given a cofibration $i$, we can form a category $\mathcal{O}(i)$ of obstruction theories for $i$ in the following way. The objects are obstruction theories for $i$ (*i.e.*, pairs $(W, \alpha)$, where $W$ is an object of $\mathcal{C}$ and $\alpha$ is a function that assigns an obstruction class to each square in a functorial way). A morphism $(W, \alpha) \to (W', \alpha')$ is a weak homotopy class in $\mathcal{C}$ from $W$ to $W'$ such that $\alpha$ is given by composition of $\alpha'$ with this weak homotopy class.

The category $\mathcal{O}(i)$ is contractible if it is non-empty. This is seen in the following way. First, note that if $(W, \alpha)$ and $(W', \alpha')$ are obstruction theories for $i$, then $(W \coprod W', \alpha \coprod \alpha')$ is also an obstruction theory for $i$. Now, fix an obstruction theory $a = (W, \alpha)$ and consider the functor $a \coprod - : \mathcal{O}(i) \to \mathcal{O}(i)$. There are natural transformations $\text{id}_{\mathcal{O}(i)} \to a \coprod - \leftarrow c_a$, where $c_a : \mathcal{O}(i) \to \mathcal{O}(i)$ is the constant functor sending everything to $a$. This shows that the identity map on the nerve of $\mathcal{O}(i)$ is null-homotopic.

In this sense, obstruction theories for $i$ are unique up to a contractible family of choices. Similarly, one can define a category of rigid obstruction theories (see Subsection 4.1), and one again finds that this category is contractible if it is non-empty.

We will show in Remark 7.3 how to construct obstruction theories in practice.

## 5. Maps that Have Obstruction Theories

Having established the basic notions of obstruction theories, we study the existence of obstruction theories for certain kinds of cofibrations.

**Proposition 5.1.** *Let $i : A \to B$ be a cofibration such that $A$ is cofibrant and $B$ is weakly contractible and fibrant. Then $i$ has a fibrant obstruction theory.*

Corollary 6.4 will show that the assumption that $B$ is fibrant is unnecessary.

*Proof.* Let $p : X \to Y$ be a fibration with fibre $F$ such that $Y$ is fibrant. Suppose given a square $i \to p$. Let $p' : X' \to B$ be the pullback fibration $X \times_Y B \to B$. The fibre of $p'$ is also $F$. Lemma 2.3 applied to the square

$$\begin{array}{ccc} X' & \xrightarrow{=} & X' \\ p' \downarrow & & \downarrow \\ B & \xrightarrow{\sim} & * \end{array}$$

implies that the map $F \to X'$ is a weak equivalence; here we use that $B$ is fibrant.

Define the obstruction for $i \to p$ to be the homotopy class of the map $A \to X'$ composed with the weak homotopy inverse of $F \to X'$. This definition is functorial.

If a lift exists for $i \to p$, then the map $A \to X'$ factors through the contractible object $B$. It follows that the obstruction is null-homotopic.

Now suppose that the obstruction is null-homotopic. We need only show that the square $i \to p'$ has a lift. Because $A$ is cofibrant and $X'$ is fibrant, we have a left null-homotopy of the map $A \to X'$. Thus, $A \to X'$ factors through a contractible object $C$. Here $C$ is the cofibre of a cofibration $A \to \mathrm{Cyl}(A)$, where $\mathrm{Cyl}(A)$ is a cylinder object for $A$. By factoring $C \to X'$ into an acyclic cofibration followed by a fibration, we may assume that $C \to X'$ is a fibration. This gives us a commutative square

$$\begin{array}{ccc} A & \longrightarrow & C \\ i \downarrow & & \downarrow \sim \\ B & \xrightarrow{=} & B \end{array}$$

in which the right vertical arrow is an acyclic fibration because both $C$ and $B$ are weakly contractible. Hence a lift exists, and the composition $B \to C \to X'$ is the desired lift. $\square$

**Proposition 5.2.** *The class of cofibrations that have an obstruction theory* (resp., *fibrant obstruction theory, cofibrant fibrant obstruction theory*) *is closed under cobase change.*

*Proof.* We prove the proposition for obstruction theories. The other cases are identical.

Let

$$\begin{array}{ccc} A & \longrightarrow & A' \\ i \downarrow & & \downarrow i' \\ B & \longrightarrow & B' \end{array}$$

be a pushout square in which $i$ (and hence $i'$ also) is a cofibration. Suppose that $i$ has an obstruction theory.

Given a map $i' \to p$ in which $p : X \to Y$ is a fibration, define the obstruction $\alpha$ to be the obstruction for the composition $i \to i' \to p$. This definition is functorial.

If $B' \to X$ is a lift for the square $i' \to p$, then the composition $B \to B' \to X$ is a lift for the square $i \to i' \to p$. On the other hand, if $B \to X$ is a lift for the square $i \to i' \to p$, then the maps $B \to X$ and $A' \to X$ induce a lift $B' \to X$ for the square $i' \to p$. Thus, a lift exists for the square $i' \to p$ if and only if a lift exists for the square $i \to i' \to p$. A lift exists for the square $i \to i' \to p$ if and only if $\alpha$ is trivial. $\square$

**Proposition 5.3.** *The class of cofibrations that have an obstruction theory* (resp., *fibrant obstruction theory, cofibrant fibrant obstruction theory*) *is closed under retract.*

*Proof.* We prove the proposition for obstruction theories. The other cases are identical.

Suppose that $i : A \to B$ is a cofibration and that it has an obstruction theory. Let $i' : A' \to B'$ be a retract of $i$. Then $i'$ is a cofibration since cofibrations are closed under retract.

Given a map $i' \to p$ in which $p : X \to Y$ is a fibration, define the obstruction $\alpha$ to be the obstruction for the composition $i \to i' \to p$. This definition is functorial.

If $B' \to X$ is a lift for the square $i' \to p$, then the composition $B \to B' \to X$ is a lift for the square $i \to i' \to p$. On the other hand, if $B \to X$ is a lift for the square $i \to i' \to p$, then the composition $B' \to B \to X$ is a lift for the square $i' \to p$. Thus, a lift exists for the square $i' \to p$ if and only if a lift exists for the square $i \to i' \to p$. A lift exists for the square $i \to i' \to p$ if and only if $\alpha$ is trivial. $\square$

## 6. Weakly Equivalent Cofibrations

The goal of this section is to show that a cofibration has an obstruction theory if and only if a weakly equivalent cofibration has an obstruction theory (Proposition 6.3). Together with the results of the previous section, this establishes the characteristic properties of the class of cofibrations that have obstruction theories.

Note that we study fibrant obstruction theories in this section. When $\mathcal{C}$ is right proper, Proposition 4.4 tells us that the same results hold for obstruction theories.

**Proposition 6.1.** *Let $i : A \to B$ and $i' : A' \to B'$ be cofibrations, and let $i \to i'$ be a weak equivalence. Also assume that $A$ and $A'$ are cofibrant or that $\mathcal{C}$ is left proper. Then $i$ has a fibrant obstruction theory if and only if $i'$ has a fibrant obstruction theory.*

*Proof.* First suppose that $i$ has a fibrant obstruction theory. Given a map $i' \to p$ in which $p : X \to Y$ is a fibration with fibrant target, define the obstruction $\alpha$ to be the obstruction for the composition $i \to i' \to p$. This definition is functorial. By Proposition 3.1, a lift exists in the square $i' \to p$ if and only if a lift exists in the square $i \to i' \to p$. A lift exists for the square $i \to i' \to p$ if and only if $\alpha$ is trivial. This finishes one implication.

Now suppose that $i'$ has a fibrant obstruction theory. By Proposition 4.5, it suffices to show that $i$ has a cofibrant fibrant obstruction theory. Let $i \to p$ be an injective cofibration such that $p : X \to Y$ is a fibration with fibrant target. Let $p'$ be the pushout $p \amalg_i i' : X \amalg_A A' \to Y \amalg_B B'$. The map $p \to p'$ is a weak equivalence because the map $i \to p$ is an injective cofibration. Here we use that the injective structure on $\text{Ar}\mathcal{C}$ is left proper or that $i$ and $i'$ are injective cofibrant so that Lemma 2.2 applies to the injective model structure on $\text{Ar}\mathcal{C}$.

The map $p'$ is not in general a fibration, so let $p''$ be an injective fibrant replacement for $p'$. Then we have a commutative diagram

$$\begin{array}{ccc} i & \longrightarrow & p \\ \sim \downarrow & & \downarrow \sim \\ i' & \longrightarrow p' \xrightarrow{\sim} & p''. \end{array}$$

Define the obstruction for $i \to p$ to be the obstruction for $i' \to p''$. This definition is functorial because the construction of $p''$ is functorial.

The obstruction vanishes if and only if $i' \to p''$ has a lift. By Proposition 3.1, $i' \to p''$ has a lift if and only if $i \to p''$ has a lift; here we use that $\mathcal{C}$ is left proper or that $A$ and $A'$ are cofibrant. By Proposition 3.2, $i \to p''$ has a lift if and only if $i \to p$ has a lift; here we use that $p$ and $p''$ have fibrant targets. □

Recall that a projective cofibrant replacement $\tilde{f} : \tilde{X} \to \tilde{Y}$ of any map $f : X \to Y$ is a cofibration between cofibrant objects together with a commuting square

$$\begin{array}{ccc} \tilde{X} & \xrightarrow{\sim} & X \\ \tilde{f} \downarrow & & \downarrow f \\ \tilde{Y} & \xrightarrow{\sim} & Y \end{array}$$

in which the horizontal arrows are acyclic fibrations.

**Corollary 6.2.** *Let $i$ be a cofibration, and suppose that $i$ has cofibrant source or that $\mathcal{C}$ is left proper. Let $\tilde{i}$ be any projective cofibrant replacement for $i$. Then $i$ has a fibrant obstruction theory if and only if $\tilde{i}$ does.*

*Proof.* There is a weak equivalence $\tilde{i} \to i$, so the result follows from Proposition 6.1. □

**Proposition 6.3.** *Let $i$ and $i'$ be weakly equivalent cofibrations. Suppose that $i$ and $i'$ have cofibrant sources or that $\mathcal{C}$ is left proper. Then $i$ has a fibrant obstruction theory if and only if $i'$ has a fibrant obstruction theory.*

The proof does not simply reduce to a repeated application of Proposition 6.1 to a zig-zag of weak equivalences because we need to know that the intermediate objects in the zig-zag in $\text{Ar}\mathcal{C}$ are cofibrations in order to apply that result.

*Proof.* Factor $i' \to *$ as
$$i' \overset{\sim}{\rightarrowtail} j \twoheadrightarrow *$$
where the first map is an acyclic projective cofibration and the second is a projective fibration. Note that $j$ is a cofibration with cofibrant source. Then $j$ is a projective fibrant replacement for $i'$, so there is an actual map $\tilde{i} \to j$ representing the weak equivalence between $i$ and $i'$, where $\tilde{i}$ is a projective cofibrant replacement for $i$.

By Corollary 6.2, $i$ has a fibrant obstruction theory if and only if $\tilde{i}$ has a fibrant obstruction theory. By Proposition 6.1, $\tilde{i}$ has a fibrant obstruction theory if and only if $j$ has a fibrant obstruction theory. By Proposition 6.1 again, $j$ has a fibrant obstruction theory if and only if $i'$ has a fibrant obstruction theory. □

**Corollary 6.4.** *Let $i : A \to B$ be a cofibration such that $B$ is weakly contractible. Suppose also that $A$ is cofibrant or that $\mathcal{C}$ is left proper. Then $i$ has a fibrant obstruction theory.*

The difference between this result and Proposition 5.1 is that we don't assume that $B$ is fibrant here.

*Proof.* Let $\tilde{i} : \tilde{A} \to \tilde{B}$ be a projective cofibrant replacement for $i$. Let $\tilde{B} \to \hat{B}$ be an acyclic cofibration from $\tilde{B}$ to a fibrant object. Then the composite cofibration $\tilde{A} \to \hat{B}$ is weakly equivalent to $i$, and it has cofibrant source and fibrant and weakly contractible target. Proposition 6.1 and Proposition 5.1 finish the argument. □

## 7. Classification of Cofibrations with Obstruction Theories

**Theorem 7.1.** *Consider the smallest class of cofibrations with cofibrant source that contains all cofibrations with cofibrant source and weakly contractible target and is closed under retract, weak equivalence, and cobase change. This class coincides with the class of cofibrations with cofibrant source that have a fibrant obstruction theory. If $\mathcal{C}$ is right proper, then this class also coincides with the class of cofibrations with cofibrant source that have an obstruction theory.*

*Proof.* The second claim follows from the first by Proposition 4.4. By Propositions 5.2, 5.3, and 6.3 and Corollary 6.4, we need only show that if $i$ is a cofibration with cofibrant source that has a fibrant obstruction theory, then $i$ is related to a cofibration with cofibrant source and weakly contractible target by a series of weak equivalences, retracts, and cobase changes.

Let $i : A \to B$ be a cofibration with a fibrant obstruction theory such that $A$ is cofibrant. We will construct the following commutative diagram:

$$\begin{array}{ccccccc}
\tilde{F} & \rightarrowtail & C\tilde{F} & & & & \\
\sim \downarrow & \searrow & \downarrow & & & & \\
F & \rightarrowtail & D & \longrightarrow & G & & \\
& \searrow & \downarrow & & \downarrow & & \\
A & \xrightarrow{=} A & \xrightarrow{\sim} A' & \rightarrowtail & C & \xrightarrow{\sim} & C' \\
i \downarrow & \tilde{i} \downarrow & i' \downarrow & & \downarrow & & \downarrow p \\
B & \xrightarrow{\sim} B^f & \xrightarrow{=} B^f & \xrightarrow{=} & B^f & \xrightarrow{=} & B^f.
\end{array}$$

Let $B \to B^f$ be a fibrant replacement for $B$, and let $\tilde{i} : A \to B^f$ be the composite, which is also a cofibration. Let $A \to A' \to B^f$ be a factorization of $\tilde{i}$ into an acyclic cofibration followed by a fibration $i'$. Let $F$ be the fibre of $i'$, and let $\tilde{F} \to F$ be a cofibrant replacement for $F$. Let $\tilde{F} \to C\tilde{F}$ be a cofibration with weakly contractible target, and let $D$ and $C$ be the pushouts as indicated above. The map $i' : A' \to B^f$ and the constant map $D \to B^f$ agree on $F$ and therefore induce a map $C \to B^f$. Factor the map $C \to B^f$ into an acyclic cofibration $C \to C'$ followed by a fibration $p : C' \to B^f$, and let $G$ be the fibre of $C' \to B^f$. The composite $D \to B^f$ is constant, so there is an induced map $D \to G$. This completes the construction of the above diagram.

The map $\tilde{i}$ has a fibrant obstruction theory by Proposition 6.3 because $i$ does. Because obstructions are functorial, the obstruction for lifting the square $\tilde{i} \to p$ is the composite of the obstruction for lifting the square $\tilde{i} \to i'$ with the map $F \to G$. Since after inverting weak equivalences the map $F \to G$ factors through the weakly contractible object $C\tilde{F}$, it is null in the homotopy category. Therefore the obstruction vanishes and a lift exists in the square $\tilde{i} \to p$. The cofibration $i$ is weakly equivalent to the cofibration $\tilde{i}$, which is a retract of the cofibration $A \to C'$ (because of the lift), which is weakly equivalent to $A' \to C$, which is a pushout of the cofibration $\tilde{F} \to C\tilde{F}$, which has cofibrant source and weakly contractible target. $\square$

**Theorem 7.2.** *Let $\mathcal{C}$ be left proper. Consider the smallest class of cofibrations that contains all cofibrations with weakly contractible target and is closed under retract, weak equivalence, and cobase change. This class coincides with the class of cofibrations that have a fibrant obstruction theory. If $\mathcal{C}$ is also right proper, then this class also coincides with the class of cofibrations that have an obstruction theory.*

*Proof.* The first claim follows from Theorem 7.1, since Corollary 6.2 allows us to assume that $i$ has cofibrant source. Use Proposition 4.4 for the second claim. $\square$

*Remark* 7.3. The smallest class of cofibrations containing cofibrations with weakly contractible target and closed under weak equivalences and cobase changes are the **principal cofibrations**. That is, up to weak homotopy, they are the maps that occur as a quotient map of some cofibre sequence. Therefore, the previous two theorems tell us that the retracts of principal cofibrations are the cofibrations that have obstruction theories.

In the case of principal cofibrations, one can inspect the proofs of Propositions 5.1, 5.2 and 6.1 and easily choose an object $W$. If $i$ is principal, then there is a cofibre

sequence
$$A \longrightarrow B \longrightarrow C$$
in which the second map is weakly homotopic to $i$. Then $A$ is one possible choice of the object $W$.

**Corollary 7.4.** *Let $B$ be a cofibrant object of $\mathcal{C}$. If the cofibration $* \to B$ has a fibrant obstruction theory, then $B$ is weakly equivalent to a retract of $\Sigma\Omega B$, where $\Sigma$ and $\Omega$ refer to the suspension and loop functors on the homotopy category (see Section 8.2).*

*Proof.* This follows from the proof of Theorem 7.1. In the notation of that proof, $\tilde{F}$ is $\Omega B$ and $C$ is $\Sigma\Omega B$. □

We shall see below that the converse is also true; retracts of suspensions always have obstruction theories.

**Corollary 7.5.** *Let $2 \leq n \leq \infty$. In the category of pointed simplicial sets, the map $* \to RP^n$ does not have an obstruction theory.*

*Proof.* The cohomology of $RP^n$ has non-trivial cup products, but the cohomology of a retract of a suspension has trivial cup products. □

Similar considerations show that the map from a point to any space with non-trivial cup products (such as tori) does not have an obstruction theory.

It is possible to prove Corollary 7.5 directly, providing an independent verification of Theorem 7.1.

One consequence of this corollary is that cofibrations that have obstruction theories are not closed under composition. The map $* \to RP^2$ is the composition of two maps that both have obstruction theories.

## 8. Applications

**8.1. Unpointed Spaces.** Our first application concerns the category of simplicial sets. When $n \geq 1$, the generating cofibration $\partial\Delta^n \to \Delta^n$ has an obstruction theory because the standard model structure on simplicial sets is right proper, because $\Delta^n$ is weakly contractible, and because $\partial\Delta^n$ is cofibrant. Strictly speaking, this obstruction theory applies to the category of pointed simplicial sets. However, one can take any unpointed lifting problem and turn it into a pointed lifting problem by choosing compatible basepoints (as long as $n \geq 1$ so that $\partial\Delta^n$ is non-empty).

We conclude that when $n \geq 1$, obstructions to lifting squares of the form
$$\begin{array}{ccc} \partial\Delta^n & \longrightarrow & X \\ \downarrow & & \downarrow p \\ \Delta^n & \longrightarrow & Y \end{array}$$
are elements of $\pi_{n-1}(F, *)$ for some $*$, where $F$ is a fibre of the fibration $p$. When $n = 0$, there is no obstruction theory for the cofibration $\emptyset \to \Delta^0$.

**Theorem 8.1.** *Let $n \geq 1$. A fibration $p : X \to Y$ of unbased simplicial sets has the right lifting property with respect to $\partial\Delta^n \to \Delta^n$ if and only if $\pi_{n-1}(F, *)$ is zero for every fibre $F$ of $p$ and every basepoint $*$ of $F$.*

*Proof.* First suppose that every homotopy group vanishes. Then lifts exist because the obstructions must be trivial.

Now suppose that $p$ has the right lifting property. Consider a fibre $F$ over an element $y$ of $Y$, and let $*$ be any basepoint of this fibre. Since $F$ is fibrant, every element of $\pi_{n-1}(F, *)$ is represented by an actual map $f : \partial \Delta^n \to F$. Take the square

$$\begin{array}{ccc} \partial\Delta^n & \longrightarrow & X \\ \downarrow & & \downarrow p \\ \Delta^n & \longrightarrow & Y \end{array}$$

in which the bottom horizontal arrow is the constant map with value $y$ and the top horizontal arrow is the composition of $f$ with the inclusion $F \to X$ of the fibre. A lift exists in the square, and this shows that $f$ is null-homotopic. □

8.2. **Suspensions of Cofibrations.** In any pointed model category, the **suspension $\Sigma A$** of any object $A$ is the cofibre of a cofibration

$$i^A : A \to CA,$$

where $CA$ is a **cone on $A$** (*i.e.*, a weakly contractible object together with a cofibration from $A$). If $\mathcal{C}$ is left proper or $A$ is cofibrant, a lifting argument combined with the dual of the proof of Lemma 2.3 shows that if $C'A$ is a *fibrant* cone on $A$, and $\Sigma'A$ is the corresponding suspension, then there is a weak equivalence $\Sigma A \to \Sigma'A$. Thus the choice of cone object does not affect the weak homotopy type of $\Sigma A$. If $f : A \to B$ is a map, and the cones are chosen so that $f$ extends to a map $CA \to CB$ (e.g. if $CB$ is fibrant, or the cones are chosen functorially), then we get an induced map $\Sigma f : \Sigma A \to \Sigma B$. In fact, $\Sigma f$ is the suspension of $f$ in the projective model structure on $\text{Ar}\mathcal{C}$, and therefore the above argument shows that any two constructions of $\Sigma f$ yield weakly equivalent maps provided that $f$ is a cofibration between cofibrant objects or that $\mathcal{C}$ is left proper.

By choosing a cone functorially, one makes $\Sigma$ into a functor, and, by the dual of Lemma 2.3, $\Sigma$ takes weak equivalences between cofibrant objects to weak equivalences. Thus the left derived functor of $\Sigma$, defined by applying $\Sigma$ to a cofibrant replacement, gives a well-defined functor on the homotopy category. If the model category is left proper, then the dual of Lemma 2.3 tells us that $\Sigma$ is homotopy invariant on all objects.

Dually, the **loop functor $\Omega$** takes an object $X$ to the fibre of a fibration

$$PX \to X,$$

where $PX$ is a **path object on $X$** (*i.e.*, a weakly contractible object with a fibration to $X$). The other statements above dualize as well.

**Theorem 8.2.** *Let $\mathcal{C}$ be a pointed model category, and let $i : A \to B$ be any cofibration. Suppose that $A$ is cofibrant or that $\mathcal{C}$ is left proper. Choose a model for $\Sigma i$ such that it is a cofibration. Then $\Sigma i$ has a fibrant obstruction theory. If $\mathcal{C}$ is right proper, then $\Sigma i$ has an obstruction theory.*

*Proof.* The second claim follows from the first claim because of Proposition 4.4. The left proper case follows from the other case by Corollary 6.2. Hence we may suppose that $A$ is cofibrant.

By two applications of the dual to [13, Prop. I.3.3], $\Sigma B$ is the homotopy cofibre of the map $Ci \to \Sigma A$, where $Ci$ is the mapping cone $CA \amalg_A B$. We may use $C(Ci) \amalg_{Ci} \Sigma A$ to compute $\Sigma B$, where $Ci \to C(Ci)$ is a cofibration with weakly contractible target. Thus $\Sigma i$ is weakly equivalent to a cobase change of the cofibration $Ci \to C(Ci)$. Now $C(Ci)$ is weakly contractible and $Ci$ is cofibrant since $A$ is, so Theorem 7.1 finishes the proof. □

*Remark* 8.3. From the proofs of Propositions 5.1 and 5.2, it follows that the obstructions for lifting $\Sigma i$ are homotopy classes from the mapping cone $Ci$. This object can be thought of as the desuspension of the cofibre of $\Sigma i$.

8.3. **Stable Model Categories.** Now we proceed to stable model categories.

A **stable model category** is a pointed model category in which the left derived suspension (see Section 8.2) induces an automorphism of the homotopy category; its inverse is then the right derived loop functor.

**Corollary 8.4.** *Let $\mathcal{C}$ be a stable model category. Every cofibration in $\mathcal{C}$ with cofibrant source has a fibrant obstruction theory. If $\mathcal{C}$ is left proper, then every cofibration in $\mathcal{C}$ has a fibrant obstruction theory. If $\mathcal{C}$ is left and right proper, then every cofibration has an obstruction theory.*

*Proof.* This follows from Theorem 8.2 and the fact that in a stable model category, every cofibration is weakly equivalent to the suspension of a cofibration. □

**Example 8.5.** Let "spectra" refer to Bousfield-Friedlander spectra [1] or symmetric spectra [9]. Both model categories are right proper. Consider the generating cofibration $F_m \partial \Delta_+^n \to F_m \Delta_+^n$ of spectra. The meaning of $F_m$ depends on which category of spectra we are considering. The spectrum $F_m \partial \Delta_+^n$ is cofibrant, but this fact is not essential because the stable model structures under consideration are left and right proper. The desuspension of the cofibre of this generating cofibration is weakly equivalent to the sphere spectrum $S^{n-m-1}$. Therefore, obstructions for lifting generating cofibrations of spectra are elements of stable homotopy groups.

**Theorem 8.6.** *Let $n \geq 1$ and $m \geq 0$. A fibration $p : X \to Y$ of spectra with fibre $Z$ has the right lifting property with respect to $F_m \partial \Delta_+^n \to F_m \Delta_+^n$ if and only if $\pi_{n-m-1} Z$ is zero.*

*Proof.* First suppose that every homotopy group vanishes. Then lifts exist because the obstructions must be trivial.

Now suppose that $p$ has the right lifting property. If $\partial \Delta^n$ and $\Delta^n$ are based at the 0th vertex, there is a pushout square

$$\begin{array}{ccc} F_m \partial \Delta_+^n & \longrightarrow & F_m \partial \Delta^n \\ \downarrow & & \downarrow \\ F_m \Delta_+^n & \longrightarrow & F_m \Delta^n, \end{array}$$

so $p$ has the right lifting property also with respect to the cofibration $F_m \partial \Delta^n \to F_m \Delta^n$.

Let $\alpha$ be any element of $\pi_{n-m-1} Z$. Since $F_m \partial \Delta^n$ is a cofibrant model for the sphere spectrum $S^{n-m-1}$, we can represent $\alpha$ by an actual map $f : F_m \partial \Delta^n \to Z$. Thus, we have a square

$$\begin{array}{ccc} F_m \partial \Delta^n & \longrightarrow & X \\ \downarrow & & \downarrow p \\ F_m \Delta^n & \longrightarrow & Y \end{array}$$

in which the top horizontal map is the composition of $f$ with the map $Z \to X$ and the bottom horizontal map is constant. This square has a lift $l$ by assumption. Since the bottom map is constant, $l$ factors through the fibre $Z$. Since the inclusion $Z \to X$ is monic, this shows that $f : F_m \partial \Delta^n \to Z$ factors through the weakly contractible object $F_m \Delta^n$. Hence $f$ is weakly null-homotopic, and $\alpha = 0$. This shows that $\pi_{n-m-1} Z = 0$. □

*Remark* 8.7. The proof of Theorem 8.6 is more complicated than the proof of Theorem 8.1. The difference arises from the fact that the simplicial set $\Delta^n$ is weakly contractible, while the spectrum $F_m \Delta^n_+$ is not.

*Remark* 8.8. Just as in Theorem 8.1, we must assume in Theorem 8.6 that $n \geq 1$. A fibration $p : X \to Y$ with fiber $Z$ has the right lifting property with respect to $F_m \partial \Delta^0_+ \to F_m \Delta^0_+$ if and only if the map $X_m \to Y_m$ of simplicial sets is surjective. This does not guarantee that $\pi_{-m-1} Z$ is zero.

Department of Mathematics, University of Western Ontario, London, Ontario N6A 5B7, Canada
*E-mail address*: `jdc@uwo.ca`

Department of Mathematics, University of Notre Dame, South Bend, IN 46556
*E-mail address*: `dwyer.1@nd.edu`

Department of Mathematics, University of Notre Dame, South Bend, IN 46556
*E-mail address*: `isaksen.1@nd.edu`